\documentclass[final,leqno,letterpaper]{article}

\def\eqref#1{(\ref{#1})}
\def\dsp{\displaystyle}
\def\Frac#1#2{\frac
{
 {\raise.6ex
 \hbox{$\displaystyle#1$}}
}
{
 {\lower.6ex
 \hbox{$\displaystyle#2$}}
 }
}

\usepackage{amsthm}

\theoremstyle{remark} 
 
\newtheorem{remark}{Remark}

\def\bigOxe{\sqcup \kern-2.3mm \sqcap}

\makeatother

\def\dsp{\displaystyle}
\def\Frac#1#2{\frac
{
 {\raise.6ex
 \hbox{$\displaystyle#1$}}
}
{
 {\lower.6ex
 \hbox{$\displaystyle#2$}}
 }
}

\def\sign{{\rm sign}}

\input epsf 

\def\CHFs#1#2#3{
{}_1F_1\left({a};{c};{z}\right)
}

\def\binom#1#2{
\renewcommand{\arraystretch}{0.9}
\left(
\begin{array}{c}
\begin{array}{c}\hskip-10pt#1\end{array}\\
\begin{array}{c}\hskip-10pt#2\end{array}
\end{array}
\hskip-10pt
\renewcommand{\arraystretch}{1.0}
\right)}

\def\erfc{{\rm erfc}}

\def\bigO{{\cal O}}

\def\wt{\widetilde}
\def\tfrac#1#2{{{\lower.6ex
\hbox{$\scriptstyle#1$}}\over 
{\raise.7ex
\hbox{$\scriptstyle#2$}}}}

\def\sign{{\rm sign}}

\def\erfc{{{\rm erfc}}}

\def\erfc{{\rm erfc}}

\def\wt{\widetilde}
\def\tfrac#1#2{{{\lower.6ex
\hbox{$\scriptstyle#1$}}\over 
{\raise.7ex
\hbox{$\scriptstyle#2$}}}}

\def\insil#1{}

 \title{Asymptotic inversion of the binomial and negative binomial cumulative distribution functions}
\author{
A. Gil\footnotemark[1]
 \and
J. Segura\footnotemark[4]
\and
N. M. Temme\footnotemark[5]
\\
}

\begin{document}

\maketitle

\renewcommand{\thefootnote}{\fnsymbol{footnote}}

\footnotetext[1]{Departamento de Matem\'atica Aplicada y CC. de la Computaci\'on.
ETSI Caminos. Universidad de Cantabria. 39005-Santander, Spain.   }
\footnotetext[4]{ Departamento de Matem\'aticas, Estadistica y 
        Computaci\'on. Universidad de Cantabria, 39005 Santander, Spain.  }
\footnotetext[5]{IAA, 1825 BD 25, Alkmaar, The Netherlands. Former address: Centrum Wiskunde \& Informatica (CWI), Science Park 123, 1098 XG Amsterdam,  The Netherlands.  }

\begin{abstract}
The computation and inversion of the binomial and negative binomial cumulative distribution functions play a key role in
many applications. In this paper, we explain how methods used for the central beta distribution function 
(described in  \cite{gil:2017:numbeta})  can be used to obtain asymptotic representations of these functions, and also for their inversion. 
The performance of the asymptotic inversion methods is illustrated with numerical examples.
\end{abstract}

Keywords: binomial cumulative distribution function; negative binomial cumulative distribution function; asymptotic representation; asymptotic  inversion methods.

AMS classification: 33B20, 41A60.

\section{Introduction}\label{sec:int}

The binomial and negative binomial distribution functions are used in many areas of science and engineering. In particular,
the generation of random binomial variables plays a key role in simulation
algorithms as, for example, the stochastic spatial modeling of chemical reactions \cite{Marquez:2007:binom}.
On the other hand,  the negative binomial distribution is, for example, widely used in genomic research to model gene expression data 
arising from RNA-sequences; see, for example, \cite{McCarthy:2012:negbin}, \cite{Li:2019:negbin}.

The binomial cumulative distribution function is defined by
\begin{equation}\label{eq:intro01}
P(n,p,x)=\displaystyle\sum_{k=0}^x \binom{n}{k}p^k(1-p)^{n-k},\quad 0\le p \le 1,
\end{equation}
with $x$ and $n$ positive integers, $x\le n$. The complementary function is 
\begin{equation}\label{eq:intro02}
Q(n,p,x)=\displaystyle\sum_{k=x+1}^n \binom{n}{k}p^k(1-p)^{n-k}=1-P(n,p,x).
\end{equation}

The negative binomial cumulative distribution function (also called Pascal distribution) can be given by

 \begin{equation}\label{eq:intro03}
P^{NB}(r,p,x)=\displaystyle\sum_{k=0}^x \binom{k+r-1}{r-1}p^r(1-p)^{k},\quad 0\le p \le 1,
\end{equation}
with $x$ and $r$ positive integers. The complementary function satisfies $Q^{NB}(n,p,x)=1-P^{NB}(r,p,x)$.
The definition of the negative binomial distribution can be extended to the case where the parameter 
$r$ takes positive real values. In this case, the distribution is called Polya distribution.
 
These functions are particular cases of the cumulative central beta distribution. 
This distribution function (also known as the incomplete beta function) is defined by 
\begin{equation}\label{eq:intro04}
I_y(a,b)=\Frac{1}{B(a,b)}\displaystyle\int_0^y t^{a-1}(1-t)^{b-1}\,dt,
\end{equation}
where we assume that $a$ and $b$ are real positive parameters and $0\le y \le 1$.
$B(a,b)$ is the Beta function

\begin{equation}\label{eq:intro05}
B(a,b)=\Frac{\Gamma(a)\Gamma(b)}{\Gamma(a+b)}.
\end{equation}

The relation between the binomial and the central beta distribution functions is the following
 
\begin{equation}\label{eq:intro06}
P(n,p,x)=I_{1-p}(n-x,x+1),\quad Q(n,p,x)=I_p(x+1,n-x).
\end{equation}

In order to avoid loss of significant digits by cancellation, it is always convenient to compute the smallest of 
the two functions ( $P(n,p,x)$ or
 $Q(n,p,x)$) . For this, one can use the transition point for the function $I_x(p,q)$,
which is given by $x_t\approx p/(p+q)$. In the case of the binomial
distribution, we will have  $p_t \approx (x+1)/(n+1)$. Then, if $p>p_t$ ($p<p_t$) it is better to evaluate
   $P(n,p,x)$ ($Q(n,p,x)$).

For the negative binomial, we have

\begin{equation}\label{eq:intro07}
P^{NB}(r,p,x)=I_{p}(r,x+1),\quad Q^{NB}(r,p,x)=I_{1-p}(x+1,r).
\end{equation}

In this case, the transition point will be given by $p_t \approx r/(r+x+1)$.  When $p<p_t$ ($p>p_t$) it is convenient to evaluate
   $P(n,p,x)$ ($Q(n,p,x)$).

In this paper, we explain that the methods used for the central beta distribution function (described in  \cite{gil:2017:numbeta})  can be used to obtain asymptotic representations of the binomial and negative binomial cumulative distribution functions, and also for inverting these functions. 

The inversion problem is, however, now slightly different: in  \cite{gil:2017:numbeta}
we considered  the problem of finding $y$ from  the equation $I_y(a,b)=\alpha$. 
In the present case,  the problem of inverting the binomial cumulative distribution function can be stated as follows: given $\alpha\in(0,1]$, $p\in(0,1)$, and $n$ (in the asymptotic problem a large positive integer), find the smallest positive integer $x$ such that
\begin{equation}\label{eq:intro08}
\alpha \le P(n,p,x)=\displaystyle\sum_{k=0}^x \binom{n}{k}p^k(1-p)^{n-k}\,.
\end{equation}
When we assume $x\in[1,n]$, we cannot take $\alpha$ smaller than the sum of the first two terms of the sum at the right-hand side. However, the sum of these two terms becomes very small when $n$ is large.  

In the finite sum definitions in \eqref{eq:intro01}, and so on, $x$ should be an integer, but in  the representations in \eqref{eq:intro06} and \eqref{eq:intro07}, $x$ may be real. In the inversion procedure we first assume that $x$ is a real parameter, and later we round $x$ to the smallest  integer larger than $x$.

We give in detail the results for the binomial cumulative distribution function and in a final section we will redefine some parameters to obtain the results for the  negative binomial cumulative distribution function.

\section{Results for the binomial distribution function}\label{sec:bindis}
In the Appendix, \S\ref{sec:incbeta}, we summarize earlier results for the incomplete beta function. We use these for the present case, where we need to change some notations.

We use the notation 
\begin{equation}\label{eq:rep01}
\nu= n+1,\quad \xi= \frac{x+1}{\nu}, \quad 1-\xi=\frac{n-x}{\nu},
\end{equation} 
and from \eqref{eq:intro07} and \eqref{eq:incbeta13} (with $a=x+1$ and $b=n-x$) it follows that
the representation of both binomial distributions $ P(n,p,x)$ and $ Q(n,p,x)$  in terms of the complementary error function is
\begin{equation}\label{eq:rep02}
\begin{array}{@{}r@{\;}c@{\;}l@{}}
P(n,p,x)&=&I_{1-p}(n-x,x+1)=\frac12\erfc\left(+\eta\sqrt{\nu/2}\right)+R_\nu(\eta),\\[8pt]
Q(n,p,x)&=&I_p(x+1,n-x)=\frac12\erfc\left(-\eta\sqrt{\nu/2}\right)-R_\nu(\eta),
\end{array}
\end{equation}
where the function $R_\nu(\eta) $ has the asymptotic expansion given in \eqref{eq:incbeta14}. The expansion can be obtained by using a recursive scheme given in  \eqref{eq:incbeta15} in terms of a  function $f(\eta)$ that arises when a change of the variable of integration is used; see \eqref{eq:incbeta02}, \eqref{eq:incbeta03} with final result in \eqref{eq:incbeta06}. In the present case we use 
\begin{equation}\label{eq:rep03}
 f(\zeta)=\frac{\lambda\zeta}{t-\xi},\quad f(\eta)=\frac{\lambda\eta}{p-\xi},\quad \lambda=\sqrt{\xi(1-\xi)},
\end{equation}
where $\zeta$ is defined in \eqref{eq:incbeta03} ($t$ is a variable of integration in \eqref{eq:incbeta02}) and the definition of $\eta$ becomes
\begin{equation}\label{eq:rep04}
-\tfrac12 \eta^2 = \xi \log \frac{p}{\xi}+(1-\xi) \log \frac{1-p}{1-\xi}, \quad \sign(\eta)=\sign(p-\xi).
\end{equation}

\begin{remark}\label{rem:rem01}
The choice of sign follows from the change of variables in \S\ref{sec:incbeta}.  We know that when  $p\downarrow0$ the binomial distributions approach the values $P(n,p,x)\to1$, $Q(n,p,x)\to 0$. From \eqref{eq:rep04} we see that the corresponding $\eta$ in the complementary error function tends to infinity when $p\downarrow0$, and when we take $\eta \to-\infty$, we have $\frac12\erfc\left(\eta\sqrt{\nu/2}\right)\to 1$, which is the wanted limit for $P(n,p,x)$. We see that this corresponds with the choice $\sign(\eta)=\sign(p-\xi)$. Similarly for $p\to1$, in which case we need positive values of $\eta$.
\end{remark}

Other representations that follow from \eqref{eq:rep02} and \eqref{eq:incbeta06} are
\begin{equation}\label{eq:rep04a}
\begin{array}{@{}r@{\;}c@{\;}l@{}}
Q(n,p,x)&=&\dsp{\frac{F_\nu(\eta)}{F_\nu(\infty)},\quad F_\nu(\eta)=\sqrt{{\frac \nu{2\pi}}}\int_{-\infty}^\eta
e^{-\frac12\nu\zeta^2}f(\zeta)\, d\zeta,}\\[8pt]
P(n,p,x)&=&\dsp{\frac{G_\nu(\eta)}{F_\nu(\infty)},\quad G_\nu(\eta)=\sqrt{{\frac \nu{2\pi}}}\int_{\eta}^{\infty}
e^{-\frac12\nu\zeta^2}f(\zeta)\, d\zeta,}
\end{array}
\end{equation}
where $f(\zeta)$ is given in \eqref{eq:rep03}.

We see here and in the representation of the incomplete beta function in \eqref{eq:incbeta06} a function $F_\nu(\infty)$, which  is defined in \eqref{eq:incbeta07}.  It has the large-$\nu$ asymptotic expansion given in \eqref{eq:incbeta07}. The first coefficients are as shown in  \eqref{eq:incbeta09}.

\subsection{Some expansions}\label{sec:binexp}

An expansion of $\eta$ in \eqref{eq:rep04} in powers of  $q=(p-\xi)/\lambda^2$  with $\lambda=\sqrt{\xi(1-\xi)}$ reads
\begin{equation}\label{eq:rep05}
\eta= q\lambda\left(1-\tfrac13(1-2\xi)q+\tfrac{1}{36}\left(7-19\xi+19\xi^2\right)q^2+\bigO\left(q^3\right)\right).\end{equation}
Limiting values (for fixed $\xi\in(0,1))$ are
\begin{equation}\label{eq:rep06}
\lim_{p\downarrow0}\eta=-\infty,\quad \lim_{p\uparrow1}\eta=+\infty. 
\end{equation}

We can also consider $\eta$ as a function of $\xi$.  
Limiting values (for fixed $p\in(0,1))$ are
\begin{equation}\label{eq:rep08}
\lim_{\xi\downarrow0}\eta=\sqrt{-2\log(1-p)},\quad \lim_{\xi\uparrow1}\eta=-\sqrt{-2\log p}. 
\end{equation}

\begin{figure}[tb]
\begin{center}
\begin{minipage}{4.5cm}
\hspace*{-0.5cm}
\epsfxsize=5cm \epsfbox{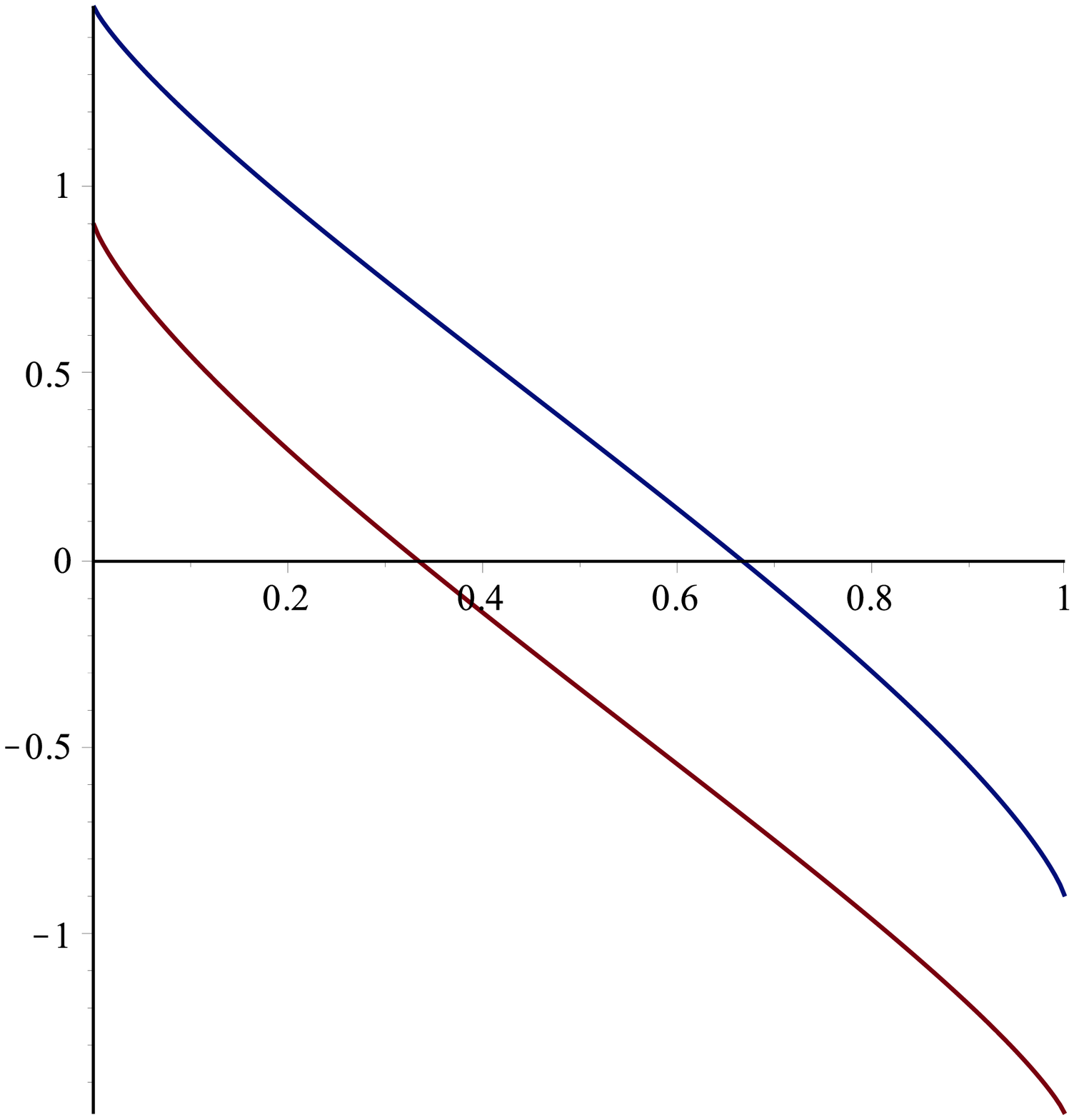}
\end{minipage}
\hspace*{2cm}
\begin{minipage}{4.5cm}
\epsfxsize=5cm \epsfbox{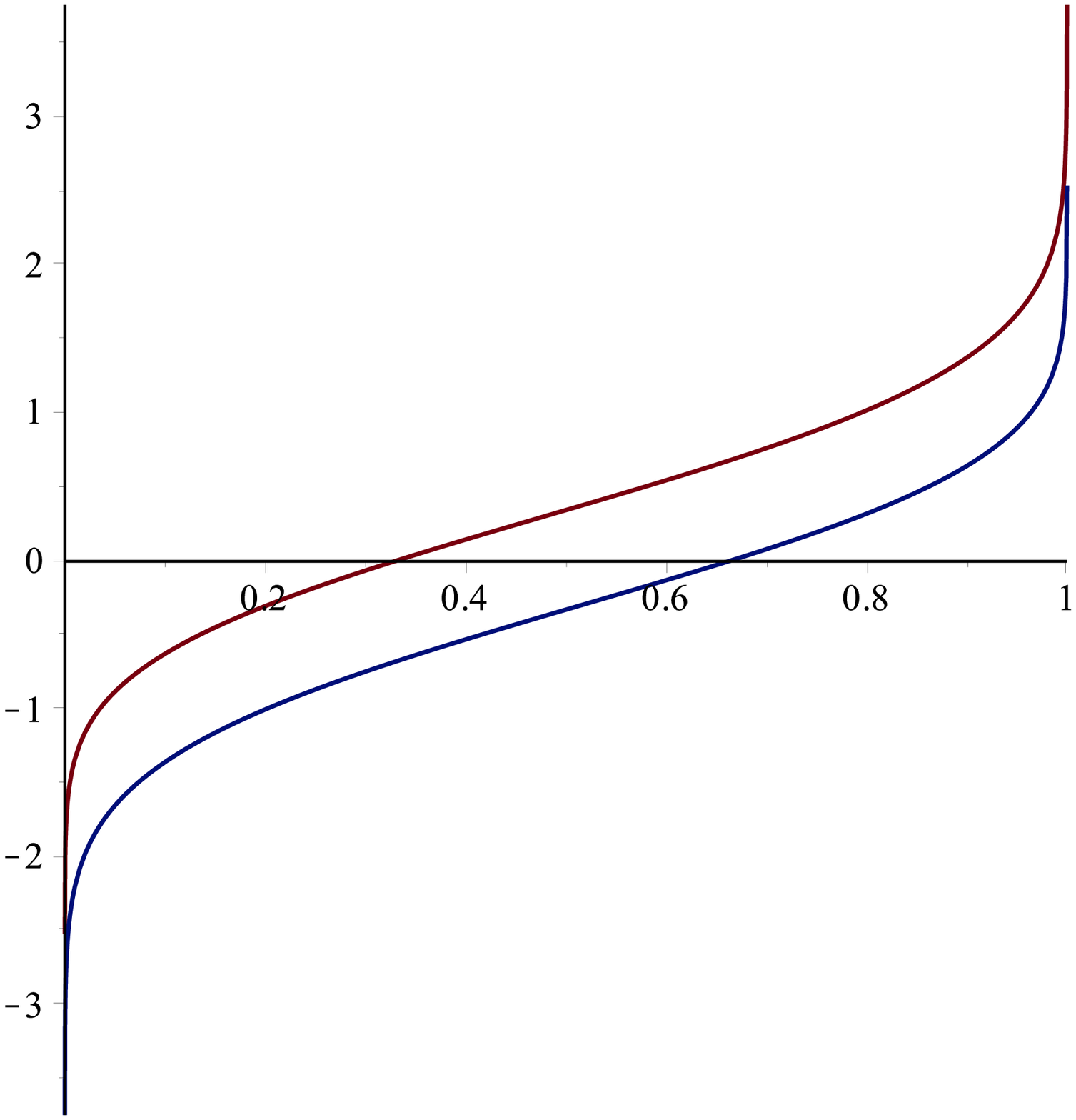}
\end{minipage}
\end{center}
\caption{{\bf Left:} The function $\eta$ defined in \eqref{eq:rep04} as a function of $\xi\in(0,1)$  for two values of $p$:
$p=1/3$ (lower curve) and $p=2/3$ (upper curve). The function $\eta$ has a zero at $\xi=p$. 
{\bf Right:} The function $\eta$ defined in \eqref{eq:rep04} as a function of $p\in(0,1)$  for two values of $\xi$:
$\xi=1/3$ (upper curve) and $\xi=2/3$ (lower curve). The function $\eta$ has a zero at $p=\xi$.}
\label{fig:fig01}
\end{figure}

In Figure~\ref{fig:fig01} (Left) we show two curves of $\eta$ as a function of $\xi$ for two values of $p$:
$p=1/3$ (upper curve) and $p=2/3$ (lower curve). The function $\eta$ has a zero at $\xi=p$. At $\xi=0$ and $\xi=1$ the values of $\eta$ follow from \eqref{eq:rep08}. In Figure~\ref{fig:fig01} (Right) we give a similar picture of $\eta$ as a function of $p$ for two values of $\xi$:
$\xi=1/3$ (lower curve) and $\xi=2/3$ upper curve). The function $\eta$ has a zero at $p=\xi$. At $p=0$ and $p=1$ we have $\eta\to\pm\infty$, see \eqref{eq:rep06}.

For the inversion procedure it is convenient to have the expansion of $\xi$ in powers of $\eta$: 
\begin{equation}\label{eq:rep09}
\xi=p-p(1-p)\sum_{k=1}^\infty a_k \wt\eta^k, \quad \wt\eta=\frac{\eta}{\sqrt{p(1-p)}}.
\end{equation}
The first coefficients are
\begin{equation}\label{eq:rep10}
\begin{array}{@{}r@{\;}c@{\;}l@{}}
a_1&=&1, \quad a_2= \frac16(2p-1),\quad a_3=\frac{1}{72}(2p^2-2p-1),\\[8pt]
a_4&=&-\frac{1}{540}(2p^3-3p^2-3p+2),\\[8pt]
a_5&=&\frac{1}{17280}(4p^4-8p^3-48p^2+52p-23).
\end{array}
\end{equation}
We also have
\begin{equation}\label{eq:rep11}
p=\xi+\lambda^2\sum_{k=1}^\infty b_k \widehat\eta^k, \quad \widehat\eta=\frac{\eta}{\lambda},\quad \lambda=\sqrt{\xi(1-\xi)},
\end{equation}
with first coefficients
\begin{equation}\label{eq:rep12}
\begin{array}{@{}r@{\;}c@{\;}l@{}}
b_1&=&1, \quad b_2= \frac13(1-2\xi),\quad b_3=\frac{1}{36}(13\xi^2-13\xi+1),\\[8pt]
b_4&=&-\frac{1}{270}(2\xi-1)(23\xi^2-23\xi-1),\\[8pt]
b_5&=&\frac{1}{4320}(313\xi^4-626\xi^3+339\xi^2-26\xi+1).
\end{array}
\end{equation}

With these coefficients we can find the coefficients of the expansion
\begin{equation}\label{eq:rep13}
f(\eta)=\frac{\lambda\eta}{p-\xi}=\sum_{k=0}^\infty c_k \widehat\eta^k, 
\end{equation}
and the first coefficients are
\begin{equation}\label{eq:rep14}
\begin{array}{@{}r@{\;}c@{\;}l@{}}
c_0&=&1, \quad c_1= \frac13(2\xi-1),\quad c_2=\frac{1}{12}(\xi^2-\xi+1),\\[8pt]
c_3&=&-\frac{1}{135}(2\xi-1)(\xi-2)(\xi+1),\quad c_4=\frac{1}{864}(\xi^2-\xi+1)^2.
\end{array}
\end{equation}

\section{Inverting the binomial distribution function using the error function}\label{sec:invbin}

We consider the inversion as described in \eqref{eq:intro08}, assuming that $\nu=n+1$ is a large parameter. The inversion procedure is based on finding $\eta$ from the equation (see \eqref{eq:rep02})
\begin{equation}\label{eq:inv01}
\tfrac12\erfc\left(\eta\sqrt{\nu/2}\right)+R_\nu(\eta)=\alpha,\quad \alpha\in(0,1),
\end{equation}
and with $\eta$ we compute $\xi$,  and then $x=\nu\xi-1$ (rounded to an integer). We consider $p$ and $n$ as fixed given quantities.

The starting point for the inversion is considering the error function in \eqref{eq:inv01} as the main term in the representation. We compute  $\eta_0$, the  solution of the reduced equation
\begin{equation}\label{eq:inv02}
\tfrac12 \erfc \left(\eta_0 \sqrt{\nu/2}\right)=\alpha.
\end{equation}

A simple and efficient algorithm for computing the inverse of
the complementary error function is included, for example, in the package described
in \cite{Gil:2015:GCH}.
Using this $\eta=\eta_0$ in \eqref{eq:rep04} we compute $\xi$, either by using the series expansion in \eqref{eq:rep09} or a numerical iteration procedure. 

\begin{remark}\label{rem:rem02}
When $\alpha$ or $1-\alpha$ is very small, the value of $\vert\eta_0\vert$ may be very large, although a large value of $\nu$ may control this. Referring to the limits  shown in \eqref{eq:rep08} for a given $p$,  we observe that if the value of $\eta_0$ satisfies $\eta_0 < -\sqrt{-2\log(1-p)}$ or $\eta_0>\sqrt{-2\log p}$, then a corresponding value of $\xi\in(0,1)$ cannot be found.
\end{remark}

Next we try to find a better approximation of $\eta$ and assume that we have an expansion of the form
\begin{equation}\label{eq:inv03}
\eta\sim\eta_0+\frac{\eta_1}{\nu}.
\end{equation}
We can find the coefficient $\eta_1$ by using a perturbation method. We have from \eqref{eq:inv02}
\begin{equation}\label{eq:inv04}
\frac{d\alpha}{d\eta_0}=-\sqrt{\frac{\nu}{2\pi}}\, e^{-\frac12\nu\eta_0^2}.
\end{equation}

To proceed,  we consider $P(n,p,x)=I_{1-p}(n-x,x+1)=\alpha$ and use the representation in \eqref{eq:rep04a}. This gives
\begin{equation}\label{eq:inv05}
\frac{d\alpha}{d\eta}=-\frac{1}{F_\nu(\infty)}\sqrt{\frac{\nu}{2\pi}}\, e^{-\frac12\nu\eta^2}f(\eta),
\end{equation}
with $f(\eta)$ given in \eqref{eq:rep03} and $\eta$ given in \eqref{eq:inv03}.

We obtain from \eqref{eq:inv04} and \eqref{eq:inv05}
\begin{equation}\label{eq:inv06}
f(\eta)\frac{d\eta}{d\eta_0}=F_\nu(\infty)e^{\frac12\nu(\eta^2-\eta_0^2)}.
\end{equation}

The coefficient $\eta_1$ in \eqref{eq:inv03} depends on $\eta_0$, and we can substitute this approximation, compare equal powers of $\nu$ and find $\eta_1$. It follows that
\begin{equation}\label{eq:inv07}
\eta_1=\frac{1}{\eta_0}\log f(\eta_0).
\end{equation}
This quantity is defined as $\eta_0\to0$ because of the expansion in \eqref{eq:incbeta10}.

For small values of $\eta_0$ (that is, when $\xi\sim p$, see \eqref{eq:rep04}), we need an expansion of $\eta_1$ in powers of $\eta_0$. We have

\begin{equation}\label{eq:inv08}
\begin{array}{@{}r@{\;}c@{\;}l@{}}
\eta_1&= &\dsp{\frac{ 1-2\xi}{ 3 \lambda}
-  \frac{5\xi^2-5\xi-1} { 36  \lambda^{2}}\eta_0-
 \frac{ (2\xi-1)(23\xi^2-23\xi-1)}
{ 1620 \lambda^{3}}  \eta_0^{2}\ -}\\ [8pt]
&& \dsp{\frac{ 31\xi^4-62\xi^3+33\xi^2-2\xi+7}{ 6480  \lambda^{4}} \eta_0^{3}+\ldots,}
\end{array}
\end{equation}
where $\lambda=\sqrt{\xi(1-\xi)}$.

\begin{remark}\label{rem:rem03}
The asymptotic estimates in this section are  uniformly valid for $\xi \in[\delta,1-\delta]$, where $\delta$ is a small fixed positive number. This corresponds with the result of the expansion of the incomplete beta function; see \eqref{eq:incbeta14}.
\end{remark}

\subsection{The algorithmic steps of the inversion procedure}\label{sec:algbin}

To summarize the algorithm 
for inverting the binomial distribution using the error function we give the following steps.

\begin{enumerate}

\item First obtain a value for $\eta$ ($\eta_0$) from \eqref{eq:inv02}.

\item With this value $\eta_0$, obtain a 
first approximation $\xi_0$ of $\xi$ from solving equation  \eqref{eq:rep04}, either by a numerical iterative procedure, or when $\eta_0$ is small  by using the expansion in \eqref{eq:rep09}.

\item Evaluate $\eta_1$ by using \eqref{eq:inv07}, where  $f(\eta_0)=\eta_0\sqrt(\xi_0(1-\xi_0))/(p-\xi_0)$; see \eqref{eq:rep03}.

\item Next compute $\eta=\eta_0+\eta_1/\nu$.

\item
With this new value of $\eta$, obtain a 
further approximation of $\xi$ by solving equation  \eqref{eq:rep04}, either by a numerical iterative procedure, or when $\eta$ is small  by using the expansion in \eqref{eq:rep09}. 

\item
Compute $x=\xi\nu-1$, and round this  to the nearest larger integer; this gives the final $x$.

\end{enumerate}

\section{Numerical examples}\label{sec:numex}
As a first example to find $x$ from $\alpha \le P(n,p,x)$, 
we take $n=50$, $p=0.4$, and $\alpha=0.51$. With $\nu=51$, we compute $\eta_0\doteq -0.0035103$ by using \eqref{eq:inv02}. This gives $\xi\doteq0.40172$ by using \eqref{eq:rep09} and $\eta_1\doteq -0.13454$ by using \eqref{eq:inv07}. Then $\eta\sim\eta_0+\eta_1/\nu\doteq-0.0061484$. The new value of $\xi$ follows from \eqref{eq:rep09}, $\xi\doteq 0.40301$. This gives $x\doteq 19.554$ and $I_{1-p}(n-x,x+1)\doteq0.510043$. 
Comparing this with $\alpha=0.51$, the absolute error is $0.000043$. Computations are done by using  Maple with Digits=16. The integer value of $x$ is 20.

When we take the same values of $\alpha$ and $p$, and $n=1500$, we find $x \doteq 599.94236$, with 
 $P(n,p,x)\doteq 0.51000026659$, an absolute error $ 2.6\times10^{-7}$. Rounding $x$ to nearest integers we find
$P(n,p,599)\doteq 0.490189$ and $P(n,p,600)\doteq 0.511212$.

\begin{figure}
\epsfxsize=13.5cm \epsfbox{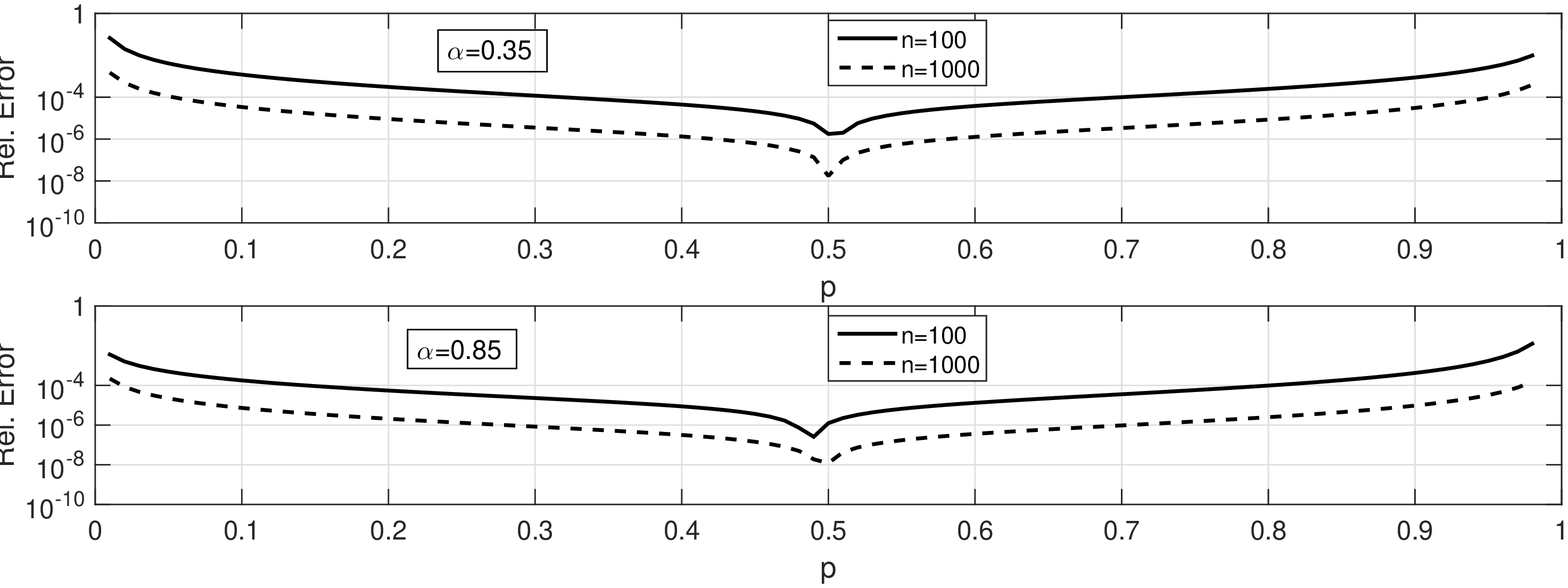}
\caption{
\label{fig:fig03} Inversion of the binomial distribution: performance of the expansion \eqref{eq:inv03}  
for $p \in (0,\,1)$ and two different values of $\alpha$ and  $n$.}
\end{figure}

A more extensive test of the performance of the expansion is considered in Figure~\ref{fig:fig03}. 
In the plots we show relative errors when the approximation  \eqref{eq:inv03}
has been considered in the inversion process for $p \in (0,\,1)$ and two different values of $\alpha$ ($\alpha=0.35,\,0.85$) and  $n$ ($n=100,\,1000$.)
As expected, a better accuracy is obtained for the larger of the two $n$-values. 

The efficiency of the computation also improves
as $n$ increases. This is not always the case in other existing algorithms for the inversion of the binomial distribution: 
for example, the CPU time in the computation of  $0.96 \le P(n,0.5,x)$ for $n=10000$ using the Matlab function {\bf binoinv}
is approximately $100$ times larger than the same computation for $n=100$. On the other hand, the algorithm implemented in R (function {\bf qbinom})
for the inversion of the binomial distribution seems to be much more efficient than the Matlab function (according
to our tests, the difference in CPU times is only a factor $2$ when computing for $n=100$ and $n=10000$) but, as before, there is not improvement
in the efficiency of the computation as $n$ increases.

\section{Results for the negative binomial distribution function}\label{sec:negbindis}
We recall the relations for the negative binomial distribution function: 
\begin{equation}\label{eq:negrep01}
P^{NB}(r,p,x)= \sum_{k=0}^x \binom{k+r-1}{r-1}p^r(1-p)^{k}=I_{p}(r,x+1),\quad 0\le p \le 1.
\end{equation}
Comparing this with the representation of $P(n,p,x)$  in \eqref{eq:intro06}, we see that we can redefine the parameters: we change $p$ into $1-p$, and write 
\begin{equation}\label{eq:negrep02}
\nu= r+x+1,\quad \xi=\frac{r}{\nu}, \quad 1-\xi=\frac{x+1}{\nu}.
\end{equation} 
The representation of the two negative binomial distributions in terms of the complementary error function is as in \eqref{eq:rep02}:
\begin{equation}\label{eq:negrep03}
\begin{array}{@{}r@{\;}c@{\;}l@{}}
P^{NB}(r,p,x)&=&I_{p}(r,x+1)=\frac12\erfc\left(-\eta\sqrt{\nu/2}\right)-R_\nu(\eta),\\[8pt]
Q^{NB}(r,p,x)&=&I_{1-p}(x+1,r)=\frac12\erfc\left(+\eta\sqrt{\nu/2}\right)+R_\nu(\eta),
\end{array}
\end{equation}
where
\begin{equation}\label{eq:negrep04}
-\tfrac12 \eta^2 = \xi \log \frac{p}{\xi}+(1-\xi) \log \frac{1-p}{1-\xi}, \quad \sign(\eta)=\sign(p-\xi).
\end{equation}
In the analysis of $P(n,p,x)$ the function $R_\nu(\eta)$ has not been used, and we refer to the Appendix to see its role in the asymptotic expansion of the incomplete beta function $I_x(a,b)$. The asymptotic  expansion of  $P^{NB}(r,p,x)$ for large $\nu$ follows from the expansion of the incomplete beta function $I_{p}(r,x+1)$.

\section{Inverting the negative binomial distribution function using the error function}\label{sec:invnegbin}
We consider the inversion problem in the form: with given positive integer $r$,  $p\in(0,1)$, and $\alpha\in(0,1)$, find the smallest integer $x$ such that
\begin{equation}\label{eq:invneg01}
\alpha \le P^{NB}(r,p,x).
\end{equation}
In particular, we assume that $r$ is large.

We use the representation in \eqref{eq:negrep03} and start with solving the equation 
\begin{equation}\label{eq:invneg02}
\tfrac12 \erfc \left(-\eta \sqrt{\nu/2}\right)=\alpha.
\end{equation}
Because the requested value of $x$ is also part of $\nu$ we have to modify the analysis for $P(n,p,x)$. We write the solution in the form
\begin{equation}\label{eq:invneg03}
-\eta \sqrt{\nu/2}=z, \quad z={\rm inverse\ erfc}(2\alpha), \quad \eta=-z\sqrt{2/\nu}=-z\sqrt{2\xi/r},
\end{equation}
because $\nu=r/\xi$. To find the corresponding $\xi$ from equation \eqref{eq:negrep04},  we write this equation in the form 
\begin{equation}\label{eq:invneg04}
\psi(\xi)=-\tfrac12\rho^2, \quad \rho=-z\sqrt{2/r}=\eta/\sqrt{\xi},  
\end{equation}
where
\begin{equation}\label{eq:invneg05}
\psi(\xi)= -\frac{1}{2\xi}\eta^2=\frac{1-\xi}{\xi}\log \frac{1-p}{1-\xi}+ \log \frac{p}{\xi},\quad \frac{d}{d\xi}\psi(\xi)=-\frac{1}{\xi^2}\log \frac{1-p}{1-\xi}.
\end{equation}
The solution $\xi$ of the equation $\psi(\xi)=-\tfrac12\rho^2$ should satisfy $\sign(p-\xi)=\sign(\eta)$.

Limiting values of the function $\psi(\xi)$ are
\begin{equation}\label{eq:invneg06}
\lim_{\xi\downarrow0}\psi(\xi)=-\infty, \quad \lim_{\xi\uparrow1}\psi(\xi)=\log p,
\end{equation}
and for $\eta$ we have
\begin{equation}\label{eq:invneg07}
\lim_{\xi\downarrow0}\eta=\sqrt{-2\log(1-p)}, \quad \lim_{\xi\uparrow1}\eta=-\sqrt{-2\log p}.
\end{equation}

So, when $\alpha <\frac12$, that is, the solution should satisfy $ p<\xi$, we can always find a solution of the equation $\psi(\xi)=-\frac12\rho^2$ for $\xi\in (0,p)$. When  $\frac12 < \alpha <1$, there is a solution for $\xi\in(p,1)$ when 
$\log p < -\frac12\rho^2$. For large values of $r$ this may be satisfied, if not we cannot use the error function  equation in 
\eqref{eq:invneg02} to find a value of $\xi$. For $p\to1$, we have $P^{NB}(r,p,x)\to1$, and the interval $(\log p,0)$ becomes very small.

For small values of $\rho$, the solution of the equation in \eqref{eq:invneg04} can be expanded in the form
\begin{equation}\label{eq:invneg08}
\xi=p-p(1-p)\sum_{k=1}^\infty r_k\wt\rho^k,\quad \wt\rho=\frac{\rho}{\sqrt{1-p}},
\end{equation}
and the first coefficients are
\begin{equation}\label{eq:invneg09}
\begin{array}{@{}r@{\;}c@{\;}l@{}}
r_1&=&1,\quad \dsp{r_2=\tfrac{1}{6}(5p-4)},\quad 
\dsp{r_3=\tfrac{1}{72}\left(47p^2-74p+26\right),}\\[8pt]
r_4&=&\dsp{\tfrac{1}{540}\left(268p^3-627p^2+453p-92\right),}\\[8pt]
r_5&=&\dsp{\tfrac{1}{17280}\left(6409p^4-19868p^3+21792p^2-9608p+1252\right).}
\end{array}
\end{equation}

We also have
\begin{equation}\label{eq:invneg10}
p=\xi+\xi(1-\xi)\sum_{k=1}^\infty s_k\widehat\rho^k,\quad \widehat\rho=\frac{\rho}{\sqrt{1-\xi}},
\end{equation}
and the first coefficients are
\begin{equation}\label{eq:invneg11}
\begin{array}{@{}r@{\;}c@{\;}l@{}}
s_1&=&1, \quad s_2= \frac13(1-2\xi),\quad s_3=\frac{1}{36}(13\xi^2-13\xi+1),\\[8pt]
s_4&=&\frac{1}{270}(1-2\xi)(23\xi^2-23\xi+1),\\[8pt]
s_5&=&\frac{1}{4320}(313\xi^4-626\xi^3+339\xi^2-26\xi+1).
\end{array}
\end{equation}

The inversion method runs as in the case for $P(n,p,x)$ with minor modifications.
\begin{enumerate}
\item
Compute $z$ and $\rho$ from \eqref{eq:invneg03} and \eqref{eq:invneg04}.
\item
Compute $\xi$ from \eqref{eq:invneg05} by solving $\psi(\xi)=-\frac12\rho^2$ by iteration or by using expansion \eqref{eq:invneg08}
when $\xi$ is small. Call this first approximation $\xi_0$ and $x_0=r/\xi_0-r-1$. 
\item The corresponding $\eta_0$ follows from equation \eqref{eq:invneg04}: $\eta_0=\rho\sqrt{\xi_0}$.
\item
Compute 
\begin{equation}\label{eq:invneg12}
\eta_1=\frac{1}{\eta_0}\log f(\eta_0), \quad f(\eta)=\frac{\eta\sqrt{\xi_0(1-\xi_0)}}{p-\xi_0}.
\end{equation}
\item
Compute $\eta=\eta_0+\eta_1/\nu$ with $\nu=r+x_0+1$.
\item
The new value $\xi$ follows from the expansion given in \eqref{eq:invneg08} 
when $\xi$ is small (or by solving $\psi(\xi)=-\frac12\rho^2$ by iteration), with $\rho=\eta/\sqrt{\xi_0}$.
\item Finally, $x=r/\xi-r-1$, rounded to the integer just larger than this value.
\end{enumerate}

As an example to find the smallest integer $x$ from $\alpha \le P^{NB}(r,p,x)$, 
we take $r=50$, $p=0.4$, and $\alpha=0.51$. The value $z$ of  \eqref{eq:invneg03} is $z\doteq -0.0177264$ and $\rho\doteq0.00354528$. Using \eqref{eq:invneg08} we obtain $\xi_0\doteq 0.398903$. Then (see \eqref{eq:invneg04})
$\eta_0=\rho\sqrt{\xi_0}\doteq0.00223916$, and \eqref{eq:invneg12} gives $\eta_1\doteq-0.137068$. With $x_0=r/\xi_0-r-1\doteq74.34369$ and $\nu\doteq125.344$. The approximation of $\eta=\eta_0+\eta_1/\nu$ becomes $\eta\doteq0.001145617$, and $\rho=\eta/\sqrt{\xi_0}\doteq0.00181387$. The corresponding $\xi$ follows from the expansion in \eqref{eq:invneg08}, which gives $\xi\doteq0.399438$, and finally $x=r/\xi-r-1\doteq74.1757$. When we compute $P^{NB}(r,p,x)$  with these values we obtain $P^{NB}(r,p,x)\doteq0.509992$. Comparing this with $\alpha$, we see an absolute error $0.79\times10^{-5}$. Computations are done by using  Maple with Digits=16.

When we take the same values of $\alpha$ and $p$, and $r=1500$, we find $x \doteq2250.71$, with 
 $P^{NB}(r,p,x)\doteq 0.50999995$, an absolute error $0.48\times10^{-7}$.

\begin{figure}
\epsfxsize=13.5cm \epsfbox{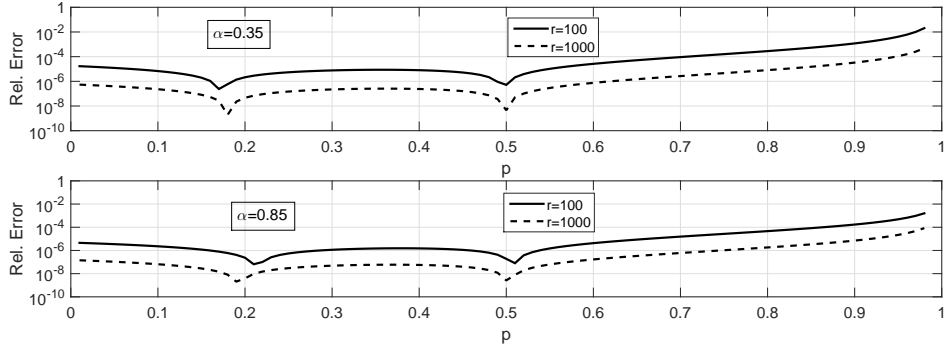}
\caption{
\label{fig:fig04} Inversion of the negative binomial distribution: performance of the expansion \eqref{eq:inv03}  
for $p \in (0,\,1)$ and two different values of $\alpha$ and  $r$.}
\end{figure}

 A more detailed example of the performance of the asymptotic inversion of the negative binomial distribution is shown 
in Figure~\ref{fig:fig04}. In the plots we show relative errors 
(obtained comparing with the values of the incomplete beta function $I_p(r,x+1)$) 
when the approximation in \eqref{eq:inv03}
has been used in the inversion process. The results obtained for 
$p\in (0,\,1)$ and
two different values of $\alpha$ ($\alpha=0.35,\,0.85$) and  $r$ ($r=100,\,1000$)
are shown for comparison.
The expansion \eqref{eq:invneg08} has been considered in all cases to obtain the value $\xi_0$.

\section{Appendix. Summary of the asymptotic results for the incomplete beta function}\label{sec:incbeta}
We collect results from  \cite{gil:2017:numbeta}, \cite{Temme:1992:AIB}, \cite[\S38.4]{Temme:2015:AMI}, with a slightly different notation. We write 
\begin{equation}\label{eq:incbeta01}
\nu=a+b,\quad \xi=\frac{a}{\nu},\quad 
b=\nu(1-\xi).
\end{equation}
Then \eqref{eq:intro04} can be written as
\begin{equation}\label{eq:incbeta02}
I_x(a,b)=\frac1{B(a,b)}\int_0^x e^{\nu\left(\xi\log
t+(1-\xi)\log(1-t)\right)}\frac{dt}{t(1-t)}.
\end{equation}
We consider $\nu$ as a large parameter, and $\xi$ bounded 
away from 0 and $1$.
The maximum of the exponential function occurs at 
$t=\xi$. We use the
transformation 
\begin{equation}\label{eq:incbeta03}
-\tfrac12\zeta^2=\xi\log\frac t{\xi}+(1-\xi)\log\frac{1-t}{1-\xi},
\end{equation}
where the sign of $\zeta$ equals the sign of $t-\xi$. 
The same transformation holds 
for
$x\mapsto\eta$ if $t$ and $\zeta$ are replaced by $x$ and 
$\eta$, respectively. That is,

\begin{equation}\label{eq:incbeta04}
-\tfrac12 \eta^2 =\xi \log \frac{x}{\xi}+(1-\xi) \log \frac{1-x}{1-\xi}. 
\end{equation}
When taking the square root for $\eta$ we assume that $\sign(\eta)=\sign(x-\xi)$, this means $\sign(\eta)=\sign\left(x-a/(a+b)\right)$.

Using 
\eqref{eq:incbeta03} we obtain
\begin{equation}\label{eq:incbeta05}
-\zeta\frac{d\zeta}{dt}=\frac{\xi-t}{t(1-t)},
\end{equation}
and we can write \eqref{eq:incbeta02} in the  form 
\begin{equation}\label{eq:incbeta06}
I_x(a,b)=\frac{F_\nu(\eta)}{F_\nu(\infty)},\quad F_\nu(\eta)=\sqrt{{\frac \nu{2\pi}}}\int_{-\infty}^\eta
e^{-\frac12\nu\zeta^2}f(\zeta)\, d\zeta,
\end{equation}
where
\begin{equation}\label{eq:incbeta07}
f(\zeta)=\frac{\zeta\lambda}{t-\xi},\quad
F_\nu(\infty)=\frac{\Gamma^*(a)\Gamma^*(b)}{\Gamma^*(a+b)}\sim
\sum_{k=0}^\infty \frac{F_k}{\nu^k}, \quad \lambda=\sqrt{\xi(1-\xi)}.
\end{equation}
The function 
$\Gamma^*(x)$, the slowly varying part of the Euler gamma function, is defined by
\begin{equation}\label{eq:incbeta08}
\Gamma^*(x)=\frac{\Gamma(x)}{\sqrt{2\pi/x}\, x^xe^{-x}},\quad x>0.
\end{equation}
The first coefficients $F_k$ are
\begin{equation}\label{eq:incbeta09}
\begin{array}{@{}r@{\;}c@{\;}l@{}}
F_0&=&1,\quad \dsp{F_1=\frac{1-\xi+\xi^2}{12\lambda^2}},\quad 
\dsp{F_2=\frac{(1-\xi+\xi^2)^2}{288\lambda^4},}\\[8pt]
F_3&=&\dsp{-\frac{139\xi^6-417\xi^5+402\xi^4-109\xi^3+402\xi^2-417\xi+139}{51840\lambda^6}.}
\end{array}
\end{equation}
The first coefficients of the Taylor expansion
\begin{equation}\label{eq:incbeta10}
f(\zeta)=a_0+a_1\zeta+a_2\zeta^2+a_3\zeta^3+\ldots
\end{equation}
are
\begin{equation}\label{eq:incbeta11}
a_0=1,\quad a_1=\frac{2\xi-1}{3\lambda},\quad 
a_2=\frac{1-\xi+\xi^2}{12\lambda^2}.
\end{equation}

When we replace in \eqref{eq:incbeta06} the function $f(\zeta)$ by 1, the integral becomes the complementary error function defined by
\begin{equation}\label{eq:incbeta12}
\erfc\,z=\frac{2}{\sqrt{\pi}}\int_z^\infty e^{-t^2}\,dt.
\end{equation}

As explained in \cite{Temme:1982:UAE}, we can write
\begin{equation}\label{eq:incbeta13}
I_x(a,b)=\tfrac12\erfc\left(-\eta\sqrt{{\nu/2}}\right)-
R_\nu(\eta),\quad \nu=a+b,
\end{equation}
where the relation between $x$ and $\eta$ follows from 
\eqref{eq:incbeta04}, and
$R_\nu(\eta)$ has the expansion 
\begin{equation}\label{eq:incbeta14}
R_\nu(\eta)\sim\frac{1}{F_\nu(\infty)} \frac{e^{-\frac12\nu\eta^2}}{\sqrt{2\pi \nu}}\sum_{k=0}^\infty\frac{C_k(\eta)}{\nu^k},\quad \nu\to\infty,
\end{equation}
and $F_\nu(\infty)$ is defined in \eqref{eq:incbeta07}. This expansion is uniformly valid for $\xi=a/(a+b)\in[\delta,1-\delta]$, where $\delta$ is a small fixed positive number.

The coefficients $C_k(\eta)$ can be obtained from the scheme
\begin{equation}\label{eq:incbeta15}
C_k(\eta)=\frac{f_k(\eta)-f_k(0)}{\eta},\quad 
f_k(\zeta)=\frac{d}{d\zeta}\frac{f_{k-1}(\zeta)-f_{k-1}(0)}{\zeta},
\end{equation}
$k=0,1,2,\ldots$, with $f_0=f$ defined in \eqref{eq:incbeta07}.

\section*{Acknowledgments}

The authors thank the anonymous referees for their constructive comments and suggestions.
This work was supported by Ministerio de Ciencia e Innovaci\'on, Spain, 
projects MTM2015-67142-P (MINECO/FEDER, UE) and PGC2018-098279-B-I00 (MCIU/AEI/FEDER, UE). 
NMT thanks CWI, Amsterdam, for scientific support.

\bibliographystyle{plain}
\bibliography{biblio}

\end{document}